\documentstyle{article}
\begin{document}

\newtheorem{theorem}{Theorem}[section]
\newtheorem{definition}{Definition}[section]
\newtheorem{corollary}[theorem]{Corollary}
\newtheorem{lemma}[theorem]{Lemma}
\newtheorem{proposition}[theorem]{Proposition}
\newtheorem{step}[theorem]{Step}
\newtheorem{example}[theorem]{Example}
\newtheorem{remark}[theorem]{Remark}

\font\sixbb=msbm6
%\font\sevenbb=msbm7
\font\eightbb=msbm8
%\font\ninebb=msbm9
%\font\tenbb=msbm10
\font\twelvebb=msbm10 scaled 1095
%\font\thirteenbb=msbm10 scaled 1315
%\font\fourteenbb=msbm10 scaled \magstep2
%%%%%%%%%%%%
\newfam\bbfam
\textfont\bbfam=\twelvebb \scriptfont\bbfam=\eightbb
                           \scriptscriptfont\bbfam=\sixbb
%\textfont\bbfam=\thirteenbb \scriptfont\bbfam=\eightbb
%                            \scriptscriptfont\bbfam=\sixbb

\newcommand{\tr}{{\rm tr \,}}
\newcommand{\linspan}{{\rm span\,}}
\newcommand{\rank}{{\rm rank\,}}
\newcommand{\diag}{{\rm Diag\,}}
\newcommand{\Image}{{\rm Im\,}}
\newcommand{\Ker}{{\rm Ker\,}}

\def\bb{\fam\bbfam\twelvebb}
\newcommand{\enp}{\begin{flushright} $\Box$ \end{flushright}}
\def\cD{{\mathcal{D}}}
\def\N{\bb N}

\title{Isometries of Grassmann spaces
\thanks{The first author was supported by the "Lend\" ulet" Program (LP2012-46/2012) of the Hungarian Academy of Sciences and by the Hungarian National Foundation for Scientific Research (OTKA), Grant No. K115383.}
\thanks{The second author was
supported by a grant from ARRS, Slovenia.}}
\author{Gy\" orgy P\' al Geh\' er \footnote{Bolyai Institute, University of Szeged, H-6720 Szeged, Aradi v\'ertan\'uk tere 1, Hungary; MTA-DE "Lend\"ulet" Functional Analysis Research Group, Institute of Mathematics University of Debrecen H-4010 Debrecen, P.O.~Box 12, Hungary. gehergy@math.u-szeged.hu, gehergyuri@gmail.com}, \quad
Peter \v Semrl\footnote{Faculty of Mathematics and Physics, University of Ljubljana,
        Jadranska 19, SI-1000 Ljubljana, Slovenia. peter.semrl@fmf.uni-lj.si}
        }

\date{}
\maketitle

\begin{abstract}
Botelho, Jamison, and Moln\' ar have recently described the general form of surjective isometries of Grassmann spaces on complex Hilbert spaces under certain dimensionality assumptions. 
In this paper we provide a new approach to this problem which enables us first, to give a shorter proof and second, to remove dimensionality constraints completely. 
In one of the low dimensional cases, which was not covered by Botelho, Jamison, and Moln\' ar, an exceptional possibility occurs. 
As a byproduct, we are able to handle the real case as well. 
Furthermore, in finite dimensions we remove the surectivity assumption. 
A variety of tools is used in order to achieve our goal, such as topological, geometrical and linear algebra techniques. 
The famous two projections theorem for two finite rank projections will be re-proven using linear algebraic methods. 
A theorem of Gy\"ory and the second author on orthogonality preservers on Grassmann spaces will be strengthened as well. 
This latter result will be obtained by using Chow's fundamental theorem of geometry of Grassmannians.
\end{abstract}
\maketitle



%-------------------------------------------------------------------------------------------------------

\section{Introduction and statement of the main results}

Describing the form of (surjective) isometries of linear normed spaces is a classical and important area of functional analysis. 
A theorem of Mazur and Ulam states that every surjective isometry between two real normed spaces is automatically an affine map, i.e. a composition of a linear map and a translation by a vector. 
Therefore if the spaces are isomorphic as metric spaces, then they are also isomorphic as vector spaces. 
Another classical example of this is the famous Banach--Stone theorem which describes surjective linear isometries between Banach spaces $C(X)$ and $C(Y)$ of continuous functions over compact Hausdorff spaces $X$ and $Y$. 
An immediate consequence of this result is that the existence of a metric isomorphism between $C(X)$ and $C(Y)$ implies the topological equivalence of the underlying spaces $X$ and $Y$. 
The non-commutative extension of this theorem was provided by  Kadison, who in particular showed that surjective linear isometries between $C^*$-algebras are closely related to algebra isomorphisms. 
We refer to \cite{FlJa,FlJa2} for more results in this direction. 

Isometries of non-linear spaces are also very important in functional analysis. 
The famous Wigner's theorem, playing an important role in the probabilistic aspects of quantum mechanics, can be interpreted as a structural result for isometries of a certain non-linear space.
Let $H$ be a complex (or real) Hilbert space. 
In quantum physics the Grassmann space $P_1(H)$ of all rank-one (orthogonal) projections is used to represent the set of pure states of the quantum system, and the quantity $\tr (PQ)$ is the so-called transition probability between two pure states. 
Wigner's theorem describes those transformations of $P_1(H)$ which preserve the transition probability. 
The conclusion is that these transformations are induced by linear or conjugate-linear isometries of $H$. 
One can easily obtain the following equation: $\|P-Q\| = \sqrt{1 - \tr PQ}$ $(P,Q\in P_1(H))$, where $\|\cdot\|$ denotes the operator norm. 
The metric on $P_1(H)$ (or on any other subset of projections) which is induced by the operator norm is usually called the gap metric. 
Therefore, Wigner's theorem characterizes isometries of $P_1(H)$ with respect to the gap metric, and in fact it states that these maps are induced by isometries of the underlying space $H$. 
Let us note that in its original version, Wigner's theorem describes surjective mappings of this kind, but as was shown later in several papers, the above conclusion holds for non-surjective transformations as well. 
The gap metric was introduced and investigated by  Sz.-Nagy and independently by  Krein and  Krasnoselski under the name "aperture". 
It has a wide range of applications from pure mathematics to engineering. 
One can easily find several references demonstrating this broad applicability, among others, we list the following fields: perturbation theory of linear operators, perturbation analysis of invariant subspaces, optimization, robust control, multi-variable control, system identification and signal processing.

Since Wigner's theorem quite a lot of attention has been paid to the study of isometries of non-linear spaces. 
Here, we are interested in the description of surjective isometries on the Grassmann space $P_n(H)$ of all rank $n$ projections with respect to the gap metric $(n\in\N)$. 
In \cite{Mo},  Moln\'ar characterized (not necessarily surjective) transformations of $P_n(H)$ which preserve the complete system of the so-called principal angles. 
These transformations are implemented by an isometry of $H$. 
The notion of principal angles was first investigated by Jordan, and has a wide range of applications such as in mathematical statistics, geometry, etc. 
We recall that the sines of the non-zero principal angles are exactly the non-zero singular values of the operator $P-Q$, each of them counted twice (see e.g.  (ii) of \cite[Theorem 26]{Gal}). 
This further implies that the quantity $\|P-Q\|$ is the sine of the largest principal angle. 
Recently,  Botelho,  Jamison, and  Moln\'ar have obtained a characterization of surjective isometries of $P_n(H)$ with respect to the gap metric for complex Hilbert spaces $H$ under the dimensionality constraint
 $\dim H \geq 4n$  (\cite{BJM}). 
Their approach was to apply a non-commutative Mazur--Ulam type result on the local algebraic behaviour of surjective isometries between substructures of metric groups. 
Then they proved that such a mapping preserves orthogonality in both directions, and finally they applied a theorem of  Gy\"ory and  \v Semrl (\cite{Gyo,Sem}), 
which contains some dimensionality constraint, too.

In the present paper, our goal is to provide a completely different approach to Botelho-Jamison-Moln\' ar's generalization of Wigner's theorem. 
We will remove the dimensionality assumption, and in finite dimensions we are able to drop the surjectivity condition. 
As a byproduct, we are also able to handle the real case. 
Furthermore, an additional possibility occurs in the case when $\dim H = 2n$, which was not covered in \cite{BJM}. 

We are now ready to state our main result on isometries of the Grassmannians.

\begin{theorem}\label{main}
Let $H$ be a complex (real) Hilbert space and $n$ a positive integer, $n < \dim H$. 
Assume that a surjective map $\phi : P_n (H) \to P_n (H)$ is an isometry with respect to the gap metric. 
If $\dim H \neq 2n$, then there exists either a unitary or an antiunitary operator (orthogonal operator) $U$ on $H$ such that $\phi$ is of the following form:
\begin{equation}\label{mainregular}
\phi (P) = UPU^\ast \qquad (P \in P_n (H)).
\end{equation}
In the case when $\dim H = 2n$,  we have either  (\ref{mainregular}), or the following additional possibility occurs:
\begin{equation}\label{mainexceptional}
\phi (P) = U(I-P)U^\ast \qquad (P \in P_n (H)).
\end{equation}
Moreover, if $\dim H < \infty$, then we have the above conclusion without assuming surjectivity.
\end{theorem}

In this paper, whenever we say a projection we automatically mean an orthogonal projection. 
Let us briefly explain our approach. 
We will consider two arbitrary projections $P$ and $Q$, and we will investigate the set $M(P,Q)$ which consists of those projections whose distance to both $P$ and $Q$ is less than or equal to $\frac{1}{\sqrt{2}}$. 
It will turn out that $M(P,Q)$ is a compact manifold if and only if $P$ and $Q$ are orthogonal. 
This will imply that orthogonality is preserved in both directions by $\phi$. 
In the case when $\dim H > 2n$ the proof is completed by a straightforward application of our second main result  stated below.
In the case when $\dim H = 2n$, the orthogonality preservers can behave badly. So, in this special case another approach is needed. It is based
on a theorem of  Blunck and  Havlicek on complementarity preservers.

We note that throughout this paper whenever we consider a manifold, we always mean a topological manifold without boundary.

Now, we state our improvement of the theorem of Gy\"ory and \v Semrl.

\begin{theorem}\label{GyS}
Let $H$ be a complex (real) Hilbert space and $n$ a positive integer such that $2n < \dim H$ is satisfied. 
Assume that a surjective map $\phi : P_n (H) \to P_n (H)$ preserves orthogonality in both directions. 
Then there exists either a unitary or an antiunitary operator (orthogonal operator) $U$ on $H$ such that 
$$
\phi (P) = UPU^\ast \qquad (P \in P_n (H)).
$$
\end{theorem}

In order to prove this result, we will apply Chow's fundamental theorem of geometry of Grassmann spaces.

In the next section we will first prove several lemmas concerning the properties of the above mentioned set $M(P,Q)$. 
We also include a proof of the well-known two projections theorem in the case when the projections are both from $P_n(H)$. 
Then Section 3 is devoted to the proofs of Theorems \ref{main} and \ref{GyS}. 

%-------------------------------------------------------------------------------------------------------

\section{Preliminary results}

In the sequel we will often use matrix representation of operators. 
In all such cases the matrices and the block-matrix forms are written with respect to an orthonormal system or an orthogonal decomposition, respectively. 
By $\diag(\dots)$ we will denote a (block-)diagonal matrix. 
In our first lemma we consider the operator norm of certain two by two matrices.

\begin{lemma}\label{twobytwo}
We have
$$
\left\| \, \left[ \begin{matrix} { {1 \over 2} & {1 \over 2} \cr {1 \over 2} & \alpha \cr} \end{matrix}\right]\, \right\| \ge {1 \over \sqrt{2} } 
\quad{\rm and}\quad
\left\| \, \left[ \begin{matrix} { {1 \over 2} & -{1 \over 2} \cr -{1 \over 2} & \alpha \cr} \end{matrix}\right]\, \right\| \ge {1 \over \sqrt{2} }
\qquad \left( - {1 \over 2} \le \alpha \le {1 \over 2} \right).
$$
Furthermore, in both cases equality holds if and only if $\alpha = -{1 \over 2}$.
\end{lemma}

\noindent
{\sl Proof}. 
The proof of the two cases are almost identical, so we will only deal with the first one. 
If we set $A_\alpha = \left[ \begin{matrix} { {1 \over 2} & {1 \over 2} \cr {1 \over 2} & \alpha \cr} \end{matrix}\right]$ with $ -{1 \over 2} \le \alpha \le {1 \over 2}$, then clearly,
\begin{equation}\label{trv}
{1 \over 2} + \alpha = \tr A_\alpha \ge 0 \ \ \ {\rm and} \ \ \ {\alpha \over 2} - {1 \over 4} = \det A_\alpha \le 0.
\end{equation}
Since $A_\alpha$ is hermitian in the complex case and symmetric in the real case, we have $\| A_\alpha \| = \max \{ | t_1 | , |t_2 | \}$, where $t_1, t_2$ are the (possibly equal) eigenvalues of $A_\alpha$. 
Because of (\ref{trv}) we get
$$
f(\alpha) = 2 \| A_\alpha \| = 2 \max \{ | t_1 | , |t_2 | \} = t_1 + t_2 + |t_1-t_2| 
$$
$$
= \tr A_\alpha + \sqrt{ (\tr A_\alpha)^2 - 4 \det A_\alpha }
= \alpha + {1 \over 2} + \sqrt{ \alpha^2 - \alpha + {5 \over 4}}.
$$
Since $f(-1/2) = \sqrt{2}$, and 
$$
f' (\alpha) = 1+\frac{-1+2 \alpha }{2 \sqrt{\frac{5}{4}-\alpha +\alpha ^2}} > 0 \qquad \left(\alpha \in \left[ - {1 \over 2}, {1 \over 2} \right]\right),
$$ 
we easily complete the proof.
\enp

The set of bounded linear operators acting on $H$ is denoted by $B(H)$. 
For any two projections $P,Q \in P_n (H)$ we define the following set:
$$
M(P,Q) = \left\{ R \in P_n (H)\, : \, \| R - P \| \le {1 \over \sqrt{2}} \ \, {\rm and}\ \, \| R - Q \| \le {1 \over \sqrt{2}} \right\}.
$$
This set will play an important role. 
If $A\subset H$ is a set, then $A^\perp$ and $\linspan A$ denote the set of all vectors which are orthogonal to every element of $A$, and the (not necessarily closed) linear manifold generated by $A$, respectively. 
In what follows, we will give a useful description of the set $M(P,Q)$ when $\|P-Q\| = 1$. 

\begin{lemma}\label{oper}
Let $H$ be a real or complex Hilbert space, $n$ a positive integer, and $P,Q \in P_n (H)$ such that $\| P - Q\| = 1$. 
Then for every $R\in M(P,Q)$ there exist an orthogonal decomposition $H = H_1 \oplus H_2$ with $\dim H_1 = 2$ and an orthonormal basis $\{ e_1 , e_2 \}$ in $H_1$, such that with respect to this decomposition and this orthonormal basis the projections $P,Q,R$ have the following matrix representations:
$$
P = \left[ \begin{matrix} { \left[ \begin{matrix} { 1 & 0 \cr 0 & 0 \cr} \end{matrix}\right] & 0 \cr 0 & P_1 \cr} \end{matrix} \right], \ \ \ 
Q = \left[ \begin{matrix} { \left[ \begin{matrix} { 0 & 0 \cr 0 & 1 \cr} \end{matrix}\right] & 0 \cr 0 & Q_1 \cr} \end{matrix} \right], \ \ \ {\rm and} \ \ \
R = \left[ \begin{matrix} { \left[ \begin{matrix} { {1 \over 2} & {1 \over 2} \cr {1 \over 2} & {1 \over 2} \cr} \end{matrix}\right] & 0 \cr 0 & R_1 \cr} \end{matrix} \right],
$$
where $P_1, Q_1 , R_1 \in P_{n-1} (H_2)$.
\end{lemma}

\noindent
{\sl Proof}. 
In both, the real and the complex case, we have $1 = \| P - Q \| = \max \{ |\langle (P-Q) x,x \rangle | \, : \, x\in H, \, \| x \| = 1 \}$. 
Since $P$ and $Q$ are projections we know that
$$
0 \le \langle Px, x \rangle , \langle Qx , x \rangle \le 1
$$
holds for every unit vector $x\in H$. 
Thus, after interchanging $P$ and $Q$, if necessary, we may assume that there exists $e_1 \in H$ such that $\| e_1 \| = 1$, $\langle P e_1 , e_1 \rangle = 1$ and $\langle Qe_1 , e_1 \rangle = 0$. 
It follows that $Pe_1 = e_1$ and $Qe_1 = 0$. 
Since $P$ and $Q$ are projections, their
matrix representations with respect to the orthogonal decomposition $H = \linspan \{e_1 \} \oplus e_{1}^\perp$ are
$$
P = \left[ \begin{matrix} { 1 & 0 \cr 0 & P_2 \cr} \end{matrix} \right] \ \ \ {\rm and} \ \ \ Q = \left[ \begin{matrix} { 0 & 0 \cr 0 & Q_2 \cr} \end{matrix} \right]
$$
for some projections $P_2, Q_2$ acting on $e_{1}^\perp$. 
The corresponding matrix representation of the projection $R$ is
$$
R = \left[ \begin{matrix} { r_1 & x^\ast \cr x & R_2 \cr} \end{matrix} \right],
$$
where $r_1$ is a real number, $x$ a vector from $e_{1}^\perp$, and $R_2 \in B(e_{1}^\perp)$. 
From $R^2 = R$,
$$
{1 \over \sqrt{2}} \ge \| R-Q \| = \left\| \, \left[ \begin{matrix} { r_1 & x^\ast \cr x & R_2 - Q_2 \cr} \end{matrix} \right]\, \right\| , 
$$
and
$$
{1 \over \sqrt{2}} \ge \| R-P \| = \left\| \, \left[ \begin{matrix} { r_1 -1 & x^\ast \cr x & R_2 - P_2 \cr} \end{matrix} \right]\, \right\| 
$$
the following equalities and inequalities can be obtained:
$$
r_{1}^2 + \| x \|^2 = r_1, \ \ \ \sqrt{ r_{1}^2 + \| x \|^2} \le {1 \over \sqrt{2}}, \ \ \ {\rm and} \ \ \ 
 \sqrt{ (r_{1}-1)^2 + \| x \|^2} \le {1 \over \sqrt{2}}.
$$
These readily imply $r_1 = {1 \over 2} = \| x \|$. 
Setting $e_2 = 2x$, the matrix representations
of $P,Q,R $ with respect to the orthogonal decomposition $H = \linspan\{e_1\} \oplus \linspan\{e_2\} \oplus \{ e_1 , e_2 \}^\perp$ are
$$
P = \left[ \begin{matrix} { 1 & 0 & 0 \cr 0 & p_2 & y^\ast \cr 0 & y & P_1 } \end{matrix} \right] , \ \ \ 
Q = \left[ \begin{matrix} { 0 & 0 & 0 \cr 0 & q_2 & w^\ast \cr 0 & w & Q_1 } \end{matrix} \right] , \ \ \ {\rm and} \ \ \
R= \left[ \begin{matrix} { {1 \over 2} & {1 \over 2} & 0 \cr {1\over 2} & r_2 & z^\ast \cr 0 & z & R_1 } \end{matrix} \right] ,
$$
for some $p_2, q_2, r_2 \in [0,1]$, $y,w,z \in \{ e_1 , e_2 \}^\perp$, and some $P_1, Q_1, R_1 \in B(\{ e_1 , e_2 \}^\perp)$.
It follows from $R^2 = R$ that $r_2 = {1 \over 2}$ and $z=0$.

Let $S\in B(H)$ be the projection onto the two-dimensional subspace $\{ e_1 , e_2 \}$. 
From
$$
\| S(R-Q)S \| \le {1 \over \sqrt{2}}
$$
we conclude that
$$
\left\| \, \left[ \begin{matrix} { {1 \over 2} & {1 \over 2} \cr {1 \over 2} & {1 \over 2} - q_2 \cr} \end{matrix}\right] \, \right\| \le {1 \over \sqrt{2}},
$$
and hence, by Lemma \ref{twobytwo}, we have $q_2 = 1$. 
But then $Q \le I$ yields that $w=0$. 
Thus, $Q$ is of the desired form, and in exactly the same way we see that also $P$ is of the form as described in the conclusion of the lemma. 
\enp

Let us note that so far we have not proven that $M(P,Q)$ is non-empty. %new sentence
We only showed that if $\|P-Q\| = 1$ and $R\in M(P,Q)$, then we have the conclusion of Lemma \ref{oper}. %new sentence

In what follows, ${\cal U}_r$ will denote either the unitary group on the $r$-dimensional complex Hilbert space, or the orthogonal
group on the $r$-dimensional real Hilbert space. 
The symbols $I_r$ and $0_r$ will denote the $r$ by $r$ identity and zero matrices, respectively.

\begin{corollary}\label{lacen}
Let $P,Q \in P_n (H)$ such that $\|P-Q\| = 1$. 
Then there exists a number $1\leq r\leq n$ such that 
\begin{equation}\label{laceneq}
P = \left[ \begin{matrix} { I_r & 0 & 0 \cr 0 & 0 & 0 \cr 0 & 0 & P_1 \cr} \end{matrix}\right] \ \ \ {\rm and} \ \ \
Q = \left[ \begin{matrix} { 0 & 0 & 0 \cr 0 & I_r & 0 \cr 0 & 0 & Q_1 \cr} \end{matrix}\right]
\end{equation}
with respect to an orthogonal decomposition $H = H_1\oplus H_2 \oplus H_3$, $\dim H_1 = \dim H_2 = r$, and $P_1,Q_1 \in P_{n-r}(H_3)$, $\| P_1 - Q_1 \| < 1$ (in the case when $r=n$ we have $P_1 = Q_1 = 0$). 
Moreover, in this case $M(P,Q)$ is the set of all projections of the form 
\begin{equation}\label{MPQ}
\left[ \begin{matrix} { {1 \over 2} I_r & {1 \over 2} U & 0 \cr {1 \over 2} U^\ast & {1 \over 2} I_r & 0 \cr 0 & 0 & R_1 \cr} \end{matrix}\right],
\end{equation}
where $U \in {\cal U}_r$, $R_1 \in P_{n-r} (H_3)$, $\| R_1 - P_1 \| \le { 1 \over \sqrt{2} }$, and $\| R_1 - Q_1 \| \le { 1 \over \sqrt{2}}$.
\end{corollary}

\noindent
{\sl Proof.} We begin with verifying (\ref{laceneq}). 
After interchanging $P$ and $Q$ if necessary, we may use exactly the same arguments as at the beginning of the proof of Lemma \ref{oper} to conclude that $P$ and $Q$ are unitary (orthogonal) similar to
$$
\left[ \begin{matrix} { 1 & 0 \cr 0 & P_2 \cr} \end{matrix} \right] \ \ \ {\rm and} \ \ \  \left[ \begin{matrix} { 0 & 0 \cr 0 & Q_2 \cr} \end{matrix} \right],
$$
where $P_2$ and $Q_2$ are projections of rank $n-1$ and $n$, respectively. 
It follows that there exists a unit vector from the intersection $\Image Q_2 \cap \Ker P_2$.
In other words, $P$ and $Q$ are unitary (orthogonal) similar to
$$
\left[ \begin{matrix} { 1 & 0 & 0 \cr 0 & 0 & 0  \cr 0 & 0  & P_3 \cr} \end{matrix} \right] \ \ \ {\rm and} \ \ \  \left[ \begin{matrix} { 0 & 0 & 0 \cr 0 & 1 & 0 \cr 0& 0 & Q_3 \cr} \end{matrix} \right],
$$
where $P_3$ and $Q_3$ are projections both of rank $n-1$. 
Now, we apply the inductive approach to obtain (\ref{laceneq}). 

Next, let $R$ be of the form (\ref{MPQ}). 
An easy calculation shows that $R\in P_n(H)$. 
We observe that the upper-left two by two corners  of $P-R$ and $Q-R$ are $\frac{1}{\sqrt{2}}$-multiples of unitary (orthogonal) operators. 
Therefore $R$ is indeed in $M(P,Q)$. 

We consider a projection $R \in M(P,Q)$. 
Then, by Lemma \ref{oper} there exists a unitary (orthogonal) operator $U$ such that
$$
P = U \, \left[ \begin{matrix} { \left[ \begin{matrix} { 1 & 0 \cr 0 & 0 \cr} \end{matrix}\right] & 0 \cr 0 & P' \cr} \end{matrix} \right] \, U^\ast , \ \ \ 
Q = U \, \left[ \begin{matrix} { \left[ \begin{matrix} { 0 & 0 \cr 0 & 1 \cr} \end{matrix}\right] & 0 \cr 0 & Q' \cr} \end{matrix} \right] \, U^\ast , $$ and $$
R = U \, \left[ \begin{matrix} { \left[ \begin{matrix} { {1 \over 2} & {1 \over 2} \cr {1 \over 2} & {1 \over 2} \cr} \end{matrix}\right] & 0 \cr 0 & R' \cr} \end{matrix} \right] \, U^\ast ,
$$
where $P', Q' , R'$ are projections of rank $n-1$. 
We also have $\| R' - P' \| \le { 1 \over \sqrt{2} }$ and $\| R' - Q' \| \le { 1 \over \sqrt{2} } $. 
If $\| P ' - Q' \| < 1$ we stop here. 
Otherwise we apply Lemma \ref{oper} once again, this time for projections $P' , Q'$, and $R'\in M(P',Q')$. 
Inductively we arrive at
$$
P = V \, \diag\left( \left[ \begin{matrix} { 1 & 0 \cr 0 & 0 \cr} \end{matrix}\right], \dots, \left[ \begin{matrix} { 1 & 0 \cr 0 & 0 \cr} \end{matrix}\right], P'' \right) \, V^\ast ,
$$
$$
Q = V \, \diag\left( \left[ \begin{matrix} { 0 & 0 \cr 0 & 1 \cr} \end{matrix}\right], \dots, \left[ \begin{matrix} { 0 & 0 \cr 0 & 1 \cr} \end{matrix}\right], Q'' \right) \, V^\ast ,
$$
$$
R = V \, \diag\left( \left[ \begin{matrix} { {1 \over 2} & {1 \over 2} \cr {1 \over 2} & {1 \over 2} \cr} \end{matrix}\right], \dots, \left[ \begin{matrix} { {1 \over 2} & {1 \over 2} \cr {1 \over 2} & {1 \over 2} \cr} \end{matrix}\right], R'' \right) \, V^\ast ,
$$
for some unitary (orthogonal) operator $V$ and some projections $P'', Q'', R''$ with $\| P'' - Q'' \| < 1$ and $R''\in M(P'',Q'')$. 
Let $k$ denote the number of two by two diagonal blocks appearing in the above matrix representations of $P$, $Q$, and $R$.

After rearranging the first $2k$ elements of the orthonormal basis of $H$ we get
\begin{equation}\label{Peq}
P = \left[ \begin{matrix} { I_r & 0 & 0 \cr 0 & 0 & 0 \cr 0 & 0 & P_1 \cr} \end{matrix}\right]
= W \, \left[ \begin{matrix} { I_k & 0 & 0 \cr 0 & 0 & 0 \cr 0 & 0 & P'' \cr} \end{matrix}\right]\, W^\ast,
\end{equation}
\begin{equation}\label{Qeq}
Q = \left[ \begin{matrix} { 0 & 0 & 0 \cr 0 & I_r & 0 \cr 0 & 0 & Q_1 \cr} \end{matrix}\right]
= W \, \left[ \begin{matrix} { 0 & 0 & 0 \cr 0 & I_k & 0 \cr 0 & 0 & Q'' \cr} \end{matrix}\right]\, W^\ast,
\end{equation}
and
\begin{equation}\label{Req}
R = W \, \left[ \begin{matrix} { {1 \over 2}I_k & {1 \over 2}I_k & 0 \cr {1 \over 2}I_k & {1 \over 2}I_k & 0 \cr 0 & 0 & R'' \cr} \end{matrix}\right]\, W^\ast
\end{equation}
for some unitary (orthogonal) operator $W$. One has to be careful when reading the above three equations. Namely, the block matrix representations of operators
on the left sides of equations correspond to the direct sum decomposition $H = H_1 \oplus H_2 \oplus H_3$, while the block matrix representations of the same
operators on the right sides correspond to some possibly different direct sum decomposition of the underlying space. But already in the next step we will show that
$k=r$, and then (after changing $W$, if necessary) we may, and we will assume that the two decompositions coincide.

From (\ref{Peq}) and (\ref{Qeq}) we infer
$$
\left[ \begin{matrix} { I_r & 0 & 0 \cr 0 & -I_r & 0 \cr 0 & 0 & P_1-Q_1 \cr} \end{matrix}\right]
= W \, \left[ \begin{matrix} { I_k & 0 & 0 \cr 0 & -I_k & 0 \cr 0 & 0 & P''-Q'' \cr} \end{matrix}\right]\, W^\ast.
$$
Comparing the eigenspaces of the two sides and taking into account that the right-bottom corners have norm less than 1, we conclude that $k=r$. 
Furthermore, the representation of $W$ with respect to the decomposition $H = H_1\oplus H_2 \oplus H_3$ is
$$
W = \left[ \begin{matrix} { W_1 & 0 & 0 \cr 0 & W_2 & 0 \cr 0 & 0 & W_3 \cr} \end{matrix}\right].
$$
Finally, from (\ref{Req}), an easy calculation gives us (\ref{MPQ}) with $U = W_1W_2^* \in {\cal U}_r$ and $R_1 = W_3R''W_3^* \in P_{n-r}(H_3)$. 
This completes the proof. 
\enp

It is important to point out that we still do not know whether the phenomena $M(P,Q) = \emptyset$ can happen or not. 
Non-emptiness of $M(P,Q)$, for arbitrary $P, Q\in P_n(H)$, is a consequence of the two projections theorem, which we will prove after the following corollary. 
However, if $P$ and $Q$ are orthogonal projections, then we do know that $M(P,Q) \neq \emptyset$, which is stated below. 

Let $P,Q \in P_n (H)$ be orthogonal projections, that is, $PQ = 0$, or equivalently, $QP=0$, which is equivalent to $\Image P \perp \Image Q$. 
Then with respect to the orthogonal decomposition $H = \Image P \oplus \Image Q \oplus H_0$ the projections $P,Q$ have the following matrix
representations:
\begin{equation}\label{jezna}
P = \left[ \begin{matrix} { I_n & 0 & 0 \cr 0 & 0 & 0 \cr 0 & 0 & 0 \cr} \end{matrix}\right] \ \ \ {\rm and} \ \ \
Q = \left[ \begin{matrix} { 0 & 0 & 0 \cr 0 & I_n & 0 \cr 0 & 0 & 0 \cr} \end{matrix}\right] .
\end{equation}

\begin{corollary}\label{orthogonal}
Let $H$ be a complex or real Hilbert space, $n$ a positive integer, and $P, Q \in P_n (H)$ projections
given by (\ref{jezna}). 
Then
$$
M(P,Q) = \left\{ \left[ \begin{matrix} { {1 \over 2} I_n & {1 \over 2} U & 0 \cr {1 \over 2} U^\ast & {1 \over 2} I_n & 0 \cr 0 & 0 & 0 \cr} \end{matrix}\right] \, : \, U \in {\cal U}_n \right\}.
$$
In particular, $M(P,Q)$ is a compact manifold.
\end{corollary}

\noindent
{\sl Proof.}
The first part is a direct consequence of the previous statement, while
the second part of the conclusion follows from the well-known facts that both the orthogonal and unitary groups are compact manifolds. 
\enp

The following lemma is known as the two projections theorem (see \cite{BS,Gal}) in the special case when $P,Q \in P_n(H)$ and $\|P-Q\|<1$. 
For the sake of completeness we give a proof here. 
Of course, the case of the two projections theorem in which the latter inequality is dropped can be obtained by combining the following lemma and Corollary \ref{lacen}.

\begin{lemma}\label{rhodos}
Let $P,Q$ be projections of rank $n$ acting on a Hilbert space $H$. 
Assume that $\| P - Q \| < 1$. 
Denote the dimension of $\Image P \cap \Image Q$ by $p$ $(0 \le p \le n)$. 
Then $P$ and $Q$ are unitary (orthogonal)
similar to operators
$$
\left[ \begin{matrix} { I_p & 0 & 0 \cr 0 & E & 0 \cr 0 & 0 & 0 \cr} \end{matrix}\right] \ \ \ {\rm and} \ \ \ \left[ \begin{matrix} { I_p & 0 & 0 \cr 0 & F & 0 \cr 0 & 0 & 0 \cr} \end{matrix}\right],
$$
where $E$ and $F$ are $2(n-p) \times 2(n-p)$ matrices given by 
$$
E = \diag \left( \left[ \begin{matrix} { 1 & 0 \cr 0 & 0 \cr } \end{matrix}\right] , \dots , \left[ \begin{matrix} { 1 & 0 \cr 0 & 0 \cr } \end{matrix}\right] \right)
$$ 
and 
$$
F = \diag \left( \left[ \begin{matrix} { d_j & \sqrt{d_j ( 1- d_j)} \cr \sqrt{d_j ( 1- d_j)} & 1-d_j \cr } \end{matrix}\right] \colon 1 \leq j \leq n-p \right),
$$ 
with $0 < d_1 , \ldots, d_{n-p} <1$.
\end{lemma}

\noindent
{\sl Proof.} We set $H_1= \Image P \cap \Image Q$ and $H_2 = H_{1}^\perp$. 
With respect to the orthogonal decomposition $H = H_1 \oplus H_2$ we have
$$
P = \left[ \begin{matrix} { I_p & 0 \cr 0 & P_1 \cr} \end{matrix}\right] \ \ \ {\rm and} \ \ \ Q = \left[ \begin{matrix} { I_p & 0 \cr 0 & Q_1 \cr } \end{matrix}\right],
$$
where $P_1$ and $Q_1$ are projections of rank $n-p$. 
The subspace $H_2$ is the orthogonal sum of $H_3 = \Image P_1$ and $H_4 = \Ker P_1$. 
With respect to the decomposition $H= H_1 \oplus H_3 \oplus H_4$ the projections $P$ and $Q$ have the following matrix representations:
$$
P = \left[ \begin{matrix} { I_p & 0 & 0 \cr 0 & I_{n-p} & 0 \cr 0 & 0 & 0 \cr} \end{matrix}\right] \ \ \ {\rm and} \ \ \ Q = \left[ \begin{matrix} { I_p & 0 & 0 \cr 0 & D_1 & D_2 \cr 0 & D_{2}^\ast & D_3 \cr} \end{matrix}\right].
$$
After applying unitary (orthogonal) similarity, if necessary, we may assume with no loss of generality that $D_1$ is diagonal, $D_1 = \diag(d_1 , \ldots , d_{n-p})$.
Moreover, the rank of the submatrix $[ D_{2}^\ast \, D_3 ]$ is at most $n-p$, and therefore, the subspace $H_4$ can be decomposed into an orthogonal sum of two subspaces,
the first one being of dimension at most $n-p$, such that the corresponding matrix representations of $P$ and $Q$ are
$$
P = \left[ \begin{matrix} { I_p & 0 & 0 & 0 \cr 0 & I_{n-p} & 0 & 0 \cr 
0 & 0 & 0 & 0 \cr 0 & 0 & 0 & 0 \cr} \end{matrix}\right] \ \ \ {\rm and} \ \ \ Q = \left[ \begin{matrix} { I_p & 0 & 0 &0 \cr 0 & D_1 & E_2 & 0 \cr 0 & E_{2}^\ast & E_3 & 0 \cr 0 & 0 & 0 & 0 \cr} \end{matrix}\right].
$$
Since $Q$ is a projection, we have $0 \le D_1 \le I_{n-p}$, and because $\| I_{n-p} - D_1 \| \le \| P - Q \| < 1$ we conclude that
$0 < d_1 , \ldots, d_{n-p} \le 1$. 
Actually, we have
$0 < d_1 , \ldots, d_{n-p} < 1$, since otherwise, one of $d_1 , \ldots, d_{n-p}$, say $d_1$, would be equal to $1$, and then since $Q$ is a projection, the first row of $E_2$ and the first column of $E_{2}^\ast$
would be zero yielding that $\dim (\Image P \cap \Image Q) \ge p+1$, a contradiction.

The size of the matrix $E_2$ is $(n-p) \times k$ with $k \le n-p$. 
We claim that actually we have $k = n-p$. 
For if this was not true, it would follow
from $Q^2 = Q$ that 
\begin{equation}\label{doma}
D_{1}^2 + E_2 E_{2}^\ast = D_1,
\end{equation} 
and consequently,
the diagonal matrix $D_1 - D_{1}^2$ would not be of full rank, which is a contradiction.

We can now apply the polar decomposition $E_2 = PU$, where $U$ is unitary (orthogonal) and $P$ is positive semidefinite. 
Applying unitary (orthogonal) similarity once more, we
can assume that already $E_2$ is positive. 
But then (\ref{doma}) yields that $E_2$ is the unique positive square root of the diagonal matrix $D_1 - D_{1}^2$. 
It follows that
$$
E_2 = \diag\left( \sqrt{d_1 ( 1- d_1)}, \ldots , \sqrt{d_{n-p} ( 1- d_{n-p})} \right),
$$
and then trivially,
$$
E_3 = \diag\left( 1- d_1, \ldots , 1- d_{n-p} \right).
$$
We complete the proof by rearranging the orthonormal basis of $H$.
\enp

The general case of the two projections theorem, i.e. when we have two finite rank projections with possibly different ranks, could be obtained from the above Lemma, Corollary \ref{lacen}, and some elementary facts concerning two projections.

Let us consider the rank one projections
$$
S = \left[ \begin{matrix} { 1 & 0 \cr 0 & 0 \cr } \end{matrix}\right] \ \ \ {\rm and} \ \ \ T(d) = \left[ \begin{matrix} { d & \sqrt{d ( 1- d)} \cr \sqrt{d ( 1- d)} & 1-d \cr } \end{matrix}\right] \ \ \ (0\leq d\leq 1).
$$ 
Some easy computations give us the following equalities and inequalities:
\begin{equation}\label{STZ}
\|S-T(1/2)\| = \frac{1}{\sqrt{2}}, \ \ \ \|T(d)-T(1/2)\| = \frac{\sqrt{1-2 \sqrt{(1-d) d}}}{\sqrt{2}} \leq \frac{1}{\sqrt{2}},
\end{equation}
\begin{equation}\label{STtilZ}
\big\|S-T\big((1+\sqrt{d})/2\big)\big\| = \big\|T(d)-T\big((1+\sqrt{d})/2\big)\big\| = \frac{\sqrt{2-2 \sqrt{d}}}{2} \leq \frac{1}{\sqrt{2}}.
\end{equation}
If we combine (\ref{STZ}) (or (\ref{STtilZ})) with the two projections theorem, then we obtain that $M(P,Q)$ is indeed non-empty for every two $P,Q \in P_n(H)$. 
We point out that if $0\leq d < 1$, then there exists a positive number $\varepsilon$ such that we have 
\begin{equation}\label{MPQout}
\|S-T(1/2-\tilde\varepsilon)\| > \frac{1}{\sqrt{2}} \ \ \ {\rm and} \ \ \ \|T(d)-T(1/2-\tilde\varepsilon)\| < \frac{1}{\sqrt{2}} \ \ \ (0 < \tilde{\varepsilon} < \varepsilon).
\end{equation}
This could be verified by straightforward calculations.

Next, as a counterpart to Corollary \ref{orthogonal}, we have the following statement.

\begin{corollary}\label{radenska}
Let $H$ be a complex or real Hilbert space, $n$ a positive integer, $2n \le \dim H$, and $P,Q \in P_n (H)$. 
Assume that $\|P - Q\| = 1$ and that $P$ and $Q$ are not orthogonal. 
Then $M(P,Q)$ is not a compact manifold. Moreover, when $H$ is of infinite dimension, then $M(P,Q)$ is not even a compact set.
\end{corollary}

\noindent
{\sl Proof.} According to Corollary \ref{lacen} we may, and we will assume that $P$ and $Q$ are of the form (\ref{laceneq}). 
We need to prove that the set $M(P,Q)$ in (\ref{MPQ}) is not a compact manifold.

Using Lemma \ref{rhodos} and (\ref{STtilZ}) it is straightforward to find
$$
R = \left[ \begin{matrix} { {1 \over 2} I_r & {1 \over 2} U_1 & 0 \cr {1 \over 2} U_{1}^\ast & {1 \over 2} I_r & 0 \cr 0 & 0 & R_1 \cr} \end{matrix}\right] \in M(P,Q)
$$ 
with $\| R_1 - P_1 \| < { 1 \over \sqrt{2} }$ and $\| R_1 - Q_1 \| < { 1 \over \sqrt{2} } $.

Hence, there exists a positive real number $\varepsilon$ such that the set ${\cal U}_\varepsilon$ consisting of all projections of the form
$$
\left[ \begin{matrix} { {1 \over 2} I_r & {1 \over 2} V & 0 \cr {1 \over 2} V^\ast & {1 \over 2} I_r & 0 \cr 0 & 0 & S \cr} \end{matrix}\right],
$$ 
where $V \in {\cal U}_r$ with $\| V - U_1 \| < \varepsilon$ and $S \in P_{n-r} (H_3)$ with $\| S - R_1 \| < \varepsilon$, is an open subset of $M(P,Q)$. 
In particular, if $\dim H = \infty$, then $M(P,Q)$ is not compact at all.

We assume from now on that $H$ is finite-dimensional. 
Assume also that $M(P,Q)$ is a compact manifold. 
Having these assumptions we need to arrive at a contradiction.

In both the real and the complex cases, the topological spaces ${\cal U}_r$ and $P_{n-r} (H_3)$ are compact manifolds.
Denote their dimensions by $q_1$ and $q_2$, respectively (the exact values of $q_1$ and $q_2$ are well-known, but not important here).
We set
\begin{equation}\label{Sdef}
{\cal S}
=
\left\{
\left[ \begin{matrix} { {1 \over 2} I_r & {1 \over 2} U & 0 \cr {1 \over 2} U^\ast & {1 \over 2} I_r & 0 \cr 0 & 0 & L \cr} \end{matrix}\right]
\colon U \in {\cal U}_r, L \in P_{n-r} (H_3)
\right\}.
\end{equation}
Then $M(P,Q) \subset {\cal S}$ and ${\cal S}$ is a compact manifold of dimension $q_1 + q_2$.
Using the fact that ${\cal U}_\varepsilon$ is an open neighbourhood of $R$ in $M(P,Q)$ as well as in ${\cal S}$ we conclude that the dimension of $M(P,Q)$ is equal to $q_1 + q_2$. 

Using Lemma \ref{rhodos} and (\ref{STZ}), we can find
$$
T = \left[ \begin{matrix} { {1 \over 2} I_r & {1 \over 2} W & 0 \cr {1 \over 2} W^\ast & {1 \over 2} I_r & 0 \cr 0 & 0 & T_1 \cr} \end{matrix}\right] \in M(P,Q)
$$ 
such that $\| T_1 - P_1 \| < { 1 \over \sqrt{2} }$ and
$\| T_1 - Q_1 \| = { 1 \over \sqrt{2} } $. 
Moreover, by (\ref{MPQout}), it is possible to find in an arbitrary neighbourhood of $T_1 \in P_{n-r} (H_3)$ a projection $T_2 \in P_{n-r} (H_3)$
such that $\| T_2 - Q_1 \| > { 1 \over \sqrt{2} }$ and $\| T_2 - P_1 \| < { 1 \over \sqrt{2} } $.

Finally, the inclusion of $M(P,Q)$ into ${\cal S}$ is a continuous injective map. 
The invariance of domain theorem states that any injective and continuous map between manifolds of the same dimensions is automatically an open map. 
Applying this theorem, we conclude that $M(P,Q)$ must be an open subset of ${\cal S}$, contradicting the fact that $T \in M(P,Q)$. 
Therefore $M(P,Q)$ is not a compact manifold.
\enp

One of the main tools in the proof of our main results is Chow's fundamental theorem of geometry of Grassmann spaces \cite{Cho}.
In this paper we prefer to speak of (orthogonal) projections rather than of subspaces. 
But if we apply the obvious identification, where a subspace of dimension $n$ is identified with a projection of rank $n$ whose image is this subspace, then we arrive at the following definition of adjacency of two projections of rank $n$: projections $P,Q \in P_n (H)$ are said to be adjacent if and only if $\dim ( \Image P + \Image Q) = n+1$ which is equivalent to $\dim ( \Image P \cap \Image Q) = n-1$.

By the two projections theorem we easily conclude that $P,Q \in P_n (H)$ are adjacent if and only if they are unitary (orthogonal) similar to
operators of the following form:
\begin{equation}\label{llrr}
P = \left[ \begin{matrix} { I_{n-1} & 0 & 0 \cr 0 & 
 \left[ \begin{matrix} { 1 & 0 \cr 0 & 0 \cr } \end{matrix}\right] & 0 \cr 0 & 0 & 0 \cr} \end{matrix}\right] \end{equation} and
\begin{equation}\label{kkrr}
 Q = \left[ \begin{matrix} { I_{n-1} & 0 & 0 \cr 0 & 
\left[ \begin{matrix} { d & \sqrt{d ( 1- d)} \cr \sqrt{d ( 1- d)} & 1-d \cr } \end{matrix}\right] 
 & 0 \cr 0 & 0 & 0 \cr} \end{matrix}\right]
\end{equation}
for some real $d$, $0 \le d < 1$. 
Equivalently, we can say that $P$ and $Q$ are adjacent if and only if $\rank(P-Q) = 2$.

A semi-linear map is an additive map $A\colon H\to H$ such that there exists a field automorphism $\sigma\colon {\bb C}\to{\bb C}$ ($\sigma\colon {\bb R}\to{\bb R}$ in the real case) which satisfies $A(\lambda x) = \sigma(\lambda) x$ for every vector $x\in H$ and every number $\lambda$.
In the case of real numbers, the only automorphism is the identity, therefore every semi-linear map is linear.
In the case of complex numbers, two trivial automorphisms are the identity and the conjugation, but there are several other automorphisms.
The above mentioned Chow's theorem states that if $2n+1 \leq \dim H < \infty$, and we have a bijective map $\phi\colon P_n(H)\to P_n(H)$ which preserves adjacency in both directions, i.e. 
$$
\rank(P-Q) = 2 \; \iff \; \rank(\phi(P)-\phi(Q)) = 2 \qquad (P,Q\in P_n(H)),
$$ 
then there exists a bijective semi-linear transformation $A\colon H\to H$ such that we have 
\begin{equation}\label{regular}
\Image \phi(P) = A (\Image P) \qquad (P\in P_n(H)).
\end{equation}
If $\dim H = 2n$, then either (\ref{regular}) holds, or we have 
\begin{equation}\label{exceptional}
\Image \phi(P) = (A (\Image P))^\perp \qquad (P\in P_n(H)).
\end{equation}

For a subset ${\cal A} \subset P_n (H)$ we define the following set
$$
{\cal A}^{\top} = \{ Q \in P_n (H)\, : \, QP = 0 \ \, {\rm for}\ \, {\rm all}\ \, P \in {\cal A} \}.
$$
The last lemma of this section characterizes adjacency of two $n$-rank projections with the help of orthogonality.

\begin{lemma}\label{ledvica}
Let $n \ge 2$ and $\dim H \ge 2n +1$. 
For $P,Q \in P_n (H)$, $P\not=Q$, the following conditions are equivalent:
\begin{itemize}
\item $P$ and $Q$ are adjacent;
\item for every $R\in P_n (H) \setminus \{ P, Q \}^{\top}$ the set $(\{ R \} \cup \{P,Q\}^{\top})^{\top}$ contains at most one projection.
\end{itemize}
\end{lemma}

\noindent
{\sl Proof.} Assume first that $P$ and $Q$ are adjacent. 
Then there is no loss of generality in assuming that they are of the form (\ref{llrr}) and (\ref{kkrr}) with respect to some orthogonal decomposition $H = H_1\oplus H_2\oplus H_3$. 
It follows that
$\{ P , Q \}^{\top}$ is the set of all rank $n$ projections of the form
$$
\left[ \begin{matrix} {0_{n-1} & 0 & 0 \cr 0 & 0_2 & 0 \cr 0 & 0 & * \cr} \end{matrix} \right].
$$
Note that the size of the bottom-right corner is at least $n\times n$, and therefore, $\{ P , Q \}^{\top}$ is not empty.

Hence, if $T\in (\{P,Q\}^{\top})^{\top}$, then $\Ker T$ contains $H_3$. 
We fix an arbitrary $R \in P_n (H) \setminus \{ P,Q \}^{\top}$ and assume that $T \in (\{ R \} \cup \{P,Q\}^{\top})^{\top}$. 
Clearly, there exist a non-zero vector $x_{12} \in H_1\oplus H_2$ and another (possibly zero) one $x_3 \in H_3$ such that $x := x_{12}\oplus x_3 \in \Image R \subset \Ker T$. 
Therefore, we obtain $\Ker T = \linspan\{x_{12}\}\oplus H_3$, and conclude that either $(\{ R \} \cup \{P,Q\}^{\top})^{\top}$ is empty, or it contains only one projection, whose range is $(\linspan\{x_{12}\}\oplus H_3)^\perp$.

We consider now the case when $P$ and $Q$ are not adjacent. 
Denote $W = \Image P + \Image Q$. 
Then $\{ P, Q \}^{\top}$ is either empty and in this case it is trivial to complete the proof; 
or it is the set of all projections of rank $n$ whose matrix representation with respect to the orthogonal decomposition $H = W \oplus W^\perp$ is of the form
$$
\left[ \begin{matrix} { 0 & 0 \cr 0 & * \cr} \end{matrix}\right].
$$
Choose 
$$
R = \left[ \begin{matrix} { E_1 & 0 \cr 0 & R_1 \cr} \end{matrix}\right],
$$
with $E_1 \in P_1(W)$ and $R_1 \in P_{n-1}(W^\perp)$. 
Using the fact that $\dim W \ge n+2$, we easily conclude that the set $(\{ R \} \cup \{P,Q\}^{\top})^{\top}$ contains infinitely many rank $n$ projections.
\enp

%-------------------------------------------------------------------------------------------------------

\section{Proofs of the main results}

Now, we are in a position to verify our main results.

\bigskip

\noindent
{\sl Proof of Theorem \ref{GyS}.} The infinite dimensional case was covered in \cite{Gyo,Sem}. 
So we may assume that $2n+1 \leq \dim H < \infty$ is satisfied. 
We would like to show that $\phi$ (which is onto) is a bijective map which preserves adjacency in both directions.
Assume first that we have $\phi(P) = \phi(Q)$.
Then $R \in P_n (H)$ is orthogonal to $P$ if and only if $\phi (R)$ is orthogonal to $\phi  (P) = \phi (Q)$ which is
equivalent to the orthogonality of $R$ and $Q$. It follows easily that $P=Q$. Hence, $\phi$ is injective, and hence bijective.

Now, by Lemma \ref{ledvica}, we easily conclude that $\phi$ preserves adjacency in both directions. 
Therefore it follows from Chow's theorem that $\phi$ has the form of (\ref{regular}) with some semi-linear mapping $A\colon H\to H$. 
Let $x$ and $y$ be two non-zero orthogonal vectors in $H$. 
We consider two projections $P,Q\in P_n(H)$ such that we have $Px = x, Py = 0, Qx = 0, Qy = y$ and $PQ = 0$. 
Therefore $Ax$ and $Ay$ are also orthogonal. 
Similarly, we can conclude that if $Ax$ and $Ay$ are orthogonal, then $x$ and $y$ has to be orthogonal as well.
An easy application of Uhlhorn's theorem \cite{U} (or \cite[Corollary 1.4]{SeJFA} together with Wigner's theorem) gives that $A$ is a non-zero scalar multiple of a unitary or an antiunitary transformation (orthogonal in the real case). 
Clearly, we can choose $A$ to be unitary or antiunitary (or orthogonal in the real case). 
Finally, using the fact that $P\in P_n(H)$ implies $UPU^\ast \in P_n(H)$ with $\Image (UPU^\ast) = U (\Image P)$, our proof is done.
\enp

If $\dim H = 2n$, then we call two projections $P,Q\in P_n(H)$ complementary if $\Image P + \Image Q = H$ is fulfilled.

\bigskip

\noindent
{\sl Proof of Theorem \ref{main}.} The case when $n=1$ is the classical version of Wigner's theorem, so we will assume $n \geq 2$ throughout the proof. 

First, assume that $\dim H < \infty$ is satisfied. 
On one hand, since $P_n(H)$ is a compact manifold, its image is also compact. 
On the other hand, the domain invariance theorem ensures that $\Image \phi$ is open as well. 
Since $P_n(H)$ is connected, we conclude the bijectivity of $\phi$. 
Therefore the surjectivity assumption is indeed disposable in the finite dimensional cases.

Second, obviously $\phi$ is a homeomorphism with respect to the topology induced by the gap metric. 
We also have 
\begin{equation}\label{phiMPQ}
\phi(M(P,Q)) = M(\phi(P),\phi(Q)).
\end{equation}
If $\dim H = \infty$, then by Corollaries \ref{orthogonal} and \ref{radenska}, the projections $P$ and $Q$ are orthogonal if and only if $M(P,Q)$ is compact. 
Therefore the map $\phi$ preserves orthogonality in both directions, and the Gy\" ory--\v Semrl theorem completes the proof of this case.

Next, we assume $2n \leq \dim H < \infty$, and we show that $\phi$ preserves orthogonality in both directions. 
Let us assume the contrary, i.e. we either have $P, Q$ with $P\perp Q$ but their images are not orthogonal; or $P, Q$ are not orthogonal but $\phi(P)\perp \phi(Q)$. 
Since $\phi^{-1}$ is also a surjective isometry, it is enough to consider the second possibility. 
Then $M(\phi(P),\phi(Q))$ is a compact manifold, but $M(P,Q)$ is not, which contradicts (\ref{phiMPQ}) and the fact that both $\phi$ and $\phi^{-1}$ are continuous.

Clearly, Theorem \ref{GyS} completes the proof in the case when $2n < \dim H < \infty$. 
Next, let us suppose that $\dim H = 2n$. 
By the two projections theorem we conclude that any two elements $P, Q\in P_n(H)$ are complementary if and only if $\|(I-P)-Q\| < 1$. 
But this is equivalent to $\|\phi(I-P)-\phi(Q)\| = \|(I-\phi(P))-\phi(Q)\| < 1$, which is satisfied if and only if $\phi(P)$ and $\phi(Q)$ are complementary. 
Hence $\phi$ preserves complementarity in both directions, and a straightforward application of \cite{BH} completes the proof of this case.

It remains to consider the case when $n < d := \dim H < 2n$ case. 
Since $\|P-Q\| = \|(I-P)-(I-Q)\|$ $(P,Q\in P_n(H))$, the map $\tilde{\phi}\colon P_{d-n}(H) \to P_{d-n}(H),\, \tilde\phi(I-P) = I-\phi(P) \; (P\in P_n(H))$ is also an isometry, but on the Grassmann space $P_{d-n}(H)$. 
Because of $1 \leq 2(d-n) < d$, we obtain that $\phi$ is of the form (\ref{mainregular}).
\enp

\section*{Acknowledgement}

The authors express their thanks for the referee's remarks, especially the one about Theorem \ref{GyS}.
%-------------------------------------------------------------------------------------------------------

\end{document}